# Optimisation of Energy and Exergy of Two-Spool Turbofan Engines using Genetic Algorithms

Vin Cent Tai[1], Phen Chiak See[1], Cristinel Mares[2], and Kjetil Uhlen[1]

[1]Norwegian University of Science and Technology, Norway
[2]Brunel University, UK

July 3, 2012

Abstract. This paper presents an application of Genetic Algorithm (GA) metaheuristics to optimise the design of two-spool turbofan engines based on exergy and energy theories. The GA is used to seek the optimal values of eight parameters that define the turbofan engine. These parameters are encoded as chromosomes (a metaphor) in GA. A computer program called TurboJet-Engine Optimizer v1.0 (TJEO-1.0) was developed by the authors to perform calculations of thermodynamic properties of the engine and implement the optimisations. The caloric properties of the working fluids are obtained by using an 8th-order polynomial fluid model taken from an existing literature. The turbo-mechanical components in the turbofan engine are as well assumed to exhibit polytropic efficiencies of today's technology level. The TJEO-1.0 is integrated with Pyevolve, an open source GA optimisation framework built for used with Python programming language. The optimal design created by TJEO-1.0 was evaluated with the following criteria: (1.) energy efficiency, (2.) exergy efficiency, and (3.) combination of both of them. The results suggest that the design of turbofan engine optimised based on the third criterion is able to produce higher specific thrust (by more than 30 %), compared with the ones optimised based on the first and the second criteria.

Keywords. Aircraft design; turbofan engine; exergy efficiency; genetic algorithm; global optimisation.
ArXiv. Subject: math.OC.

## 1 Introduction

The design of aircraft engines are complex as they contain tens of thousands of components [Homaifar et al., 1994]. It begins with parametric cycle analysis [Shwin, 2010], which includes the determination of compressors' compression ratios, burner outlet temperature, bypass ratio, etc. As a preventive measure, it is important to determine the optimal value of each of these parameters in the early phase of engine development. This is because, the cost of solving a problem related to the engine design increases exponentially with the stage of production life-cycle before it (the fault) is discovered.

Many different criteria have been used to evaluate the performance of aircraft engines. For example, [Homaifar et al., 1994] used specific thrust and overall first-law efficiency as the performance measurement, [Asako et al., 2002] employed thrust specific fuel consumption (TSFC) as the performance measurement, and [Atashkari et al., 2005] measured the engine performance with TSFC, specific thrust, propulsive efficiency, and thermal efficiency.

The quantity and quality of the energy have to be taken into account for effective and efficient use of fuels [Dincer et al., 2004]. Exergy is the maximum theoretical work obtainable during a process that brings the system into equilibrium with its environment [Moran et al., 2011]. It is a very useful tool to detect and minimise the irreversible losses [Amati et al., 2006]. [Bejan and Siems, 2001] and [Riggins, 2003] pointed out the need to incorporate exergy analysis and thermodynamic optimisation in aircraft engine design. The application of exergy analysis on turbojet engine is presented in [Rosen and Etele, 2004]. [Amati et al., 2006] combined the energy and exergy methods in developing a simulation tool for scramjet design. [Tona et al., 2010] applied an exergy based analysis to evaluate the global performance of a turbofan engine and its components over a complete flight mission. All of these research pointed out exergy analysis should be included as a performance measurement for aircraft engines. This inspired the authors to incorporate the exergy analysis and the conventional thermodynamic cycle analysis in turbofan design.

This paper presents an application of exergy and energy analysis to the optimisation of two-spool turbofan engines by means of Genetic Algorithm (GA). The symbols used in this paper are given in Appendix A. After a brief introduction on the two-spool turbofan system used in this study, the concepts of energy and exergy analysis applied to turbofan performance measurements are presented. Then, the optimisation formulations are presented, followed by the results and discussions. This paper ends with conclusions.

## 2    Turbofan Engines

Figure 1 illustrates the turbofan model used in this study, alongside with the station numbers. It is a two-spool turbofan with intermediate stage bleed. Ambient air stream (station 0) enters the engine (station number 1) and separates into two streams; the bypass stream passes through the fan and exits the fan nozzle at station 7 produces cold thrust, while the core stream passes through the fan and enters the Low and High Pressure Compressors (LPC and HPC) for further compression. Cooling air streams are taken together with the bleed air in between the HPC stages. The remaining core stream is channelled to the burner (also known as combustion chamber). Fuel is injected into the core stream and is combusted in the combustion chamber. The combustion product is then cooled with a fraction of cooling air extracted previously, before it enters the High-Pressure Turbine (HPT). The gas is expanded across the HPT and the work produced is used to run the HPC and to extract the cooling air. The gas at the exit of the HPT is mixed with another fraction of cooling air for further cooling. The gas expands as it passes through the Low-Pressure Turbine (LPT) and the work produced is used to run the fan and LPC. The air then exits through the core nozzle (station 6) and produces hot thrust.

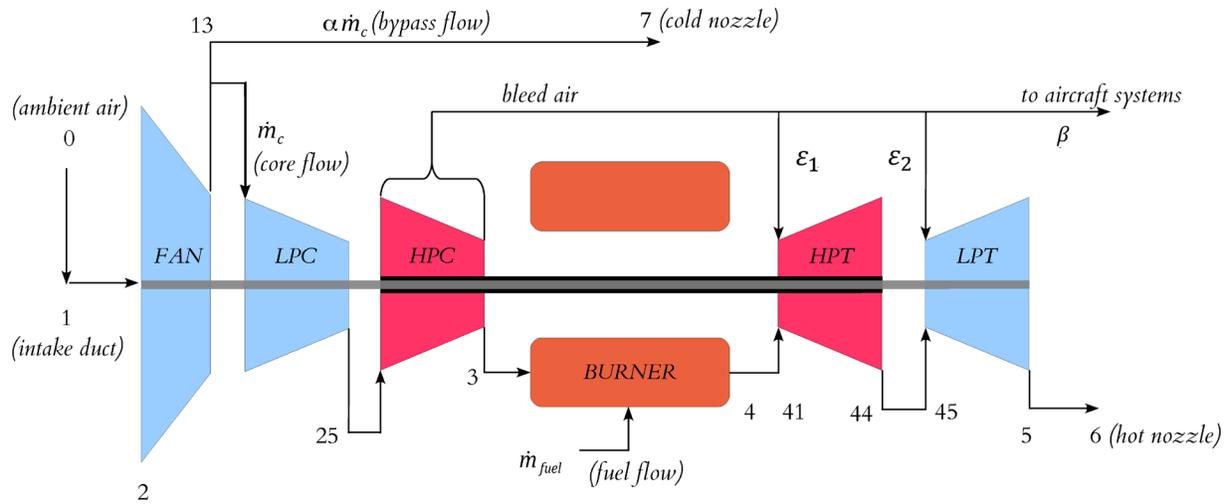

Figure 1. Schematic diagram of the investigated two-spool turbofan system

## 3 Engine performance measurements

### 3.1 Energy based performance analysis

The following assumptions were made to establish the thermodynamic analysis of the turbofan model: (1.) the working fluids (airstream and gases) are perfect; and (2.) both specific heat constant at constant pressure, $c_p$ and heat capacity ratio, $\gamma$ of the working fluids are not constant at each engine station. Further, the changes of gas constant, $R$ throughout the thermodynamic cycle due to the change in the molecular compositions of the working fluids caused by the change in fuel-air ratio, $f$ (see (4)) at each stage were also considered.

The fluid model used in this study for thermodynamic properties calculations is taken from [Walsh and Fletcher, 2004]. This model is accurate when the fluid temperature is in the range of 200-2000K [Guha, 2001]. Since the thermodynamic properties of a gas vary with temperature and its chemical composition, they can be expressed in terms of temperature and fuel-air ratio as follows:

$$c_p(T,f) = \sum_{k=0}^{8} A_k \left(\frac{T}{1000}\right)^k + \frac{f}{f+1} \sum_{k=0}^{7} B_k \left(\frac{T}{1000}\right)^k \tag{1}$$

$$h(T,f) = A_9 + \sum_{k=0}^{8} \frac{A_k}{k+1} \left(\frac{T}{1000}\right)^{k+1} + \frac{f}{f+1} \left[B_8 + \sum_{k=0}^{7} \frac{B_k}{k+1} \left(\frac{T}{1000}\right)^{k+1}\right] \tag{2}$$

$$\bar{s}(T,f) = A_{10} + A_0 \ln\left(\frac{T}{1000}\right) + \sum_{k=0}^{8} \frac{A_k}{k} \left(\frac{T}{1000}\right)^k + \frac{f}{f+1} \left[B_9 + B_0 \ln\left(\frac{T}{1000}\right) + \sum_{k=0}^{7} \frac{B_k}{k} \left(\frac{T}{1000}\right)^k\right] \tag{3}$$

where $h$ is the specific enthalpy, $\bar{s}$ is the absolute entropy at constant pressure, $A_k$ and $B_k$ are polynomial constants for dry air and Kerosene fuel respectively. Note that, when $f$ is zero, the thermodynamic properties calculated in (1)-(4) reduce to those of dry air. The values of $A_k$ and $B_k$ are listed in Table 1.

Table 1. Polynomial constants used for calculating the gas properties.

| Dry air | | Kerosene | |
| --- | --- | --- | --- |
| $A_k$ | Value | $B_k$ | Value |
| $A_0$ | 0.9923 | $B_0$ | -0.7189 |
| $A_1$ | 0.2367 | $B_1$ | 8.7475 |
| $A_2$ | -1.8524 | $B_2$ | -15.8632 |
| $A_3$ | 6.0832 | $B_3$ | 17.2541 |
| $A_4$ | -8.8940 | $B_4$ | -10.2338 |
| $A_5$ | 7.0971 | $B_5$ | 3.0818 |
| $A_6$ | -3.2347 | $B_6$ | -0.3611 |
| $A_7$ | 0.7946 | $B_7$ | 0.0039 |
| $A_8$ | -0.0819 | $B_8$ | 0.0556 |
| $A_9$ | 0.4222 | $B_9$ | -0.0016 |
| $A_{10}$ | 0.0011 | $B_{10}$ | - |

Kerosene with Lower Heating Value (LHV) equals to 43124kJ/kg was assumed. The gas constant for the aforementioned Kerosene fuel, according to [Walsh and Fletcher, 2004], changes with fuel-air ratio as follows:

$$R = 287.05 - 0.0099f + 10^{-7}f^2 \tag{4}$$

The quantity of energy of a system is evaluated with the following equation, with the assumption that neither heat nor energy is rejected from the system (as stated by the first law of thermodynamics):

$$\sum \dot{m}_{in} h_{t,in} + \sum \dot{W}_{in} = \sum \dot{m}_{out} h_{t,out} + \sum \dot{W}_{out} \tag{5}$$

The two-spool turbofan model which consists of the algorithms that compute the thermodynamic properties (based on (1)-(5)) of the working fluid at each station, are presented in Appendix B. Polytropic efficiencies of rotating components, total pressure ratios of inlet and nozzles, and the burner efficiency, as well as the shaft mechanical efficiencies were also taken into account in the analysis. Their values, based on today's technology level [Mattingly et al., 2002], are presented in Table 2.

Table 2. Efficiencies of turbofan components.

| Assumptions | | Note |
|---|---|---|
| Total pressure ratios | $\pi_d$ = 0.99 | Total pressure ratio across diffuser. |
| | $\pi_b$ = 0.96 | Total pressure ratio across burner. |
| | $\pi_{nf}$ = 0.96 | Total pressure ratio across cold nozzle. |
| | $\pi_{nc}$ = 0.96 | Total pressure ratio across core nozzle. |
| Component efficiencies | $\eta_b$ = 0.99 | Burner efficiency. |
| | $\eta_{mL}$ = 0.99 | Mechanical efficiency of low pressure spool. |
| | $\eta_{mH}$ = 0.99 | Mechanical efficiency of high pressure spool. |
| Polytropic efficiencies | $e_f$ = 0.89 | Polytropic efficiency of fan. |
| | $e_c$ = 0.90 | Polytropic efficiency of compressors. |
| | $e_t$ = 0.89 | Polytropic efficiency of turbines. |

The energy efficiency, $\eta_I$ of the engine is computed with (6):

$$\eta_I = \frac{F_6 V_6 + F_7 V_7 - \dot{m}_0 V_0}{\dot{m}_{fuel} |LHV|} \qquad (6)$$

where $F_7$ is the cold thrust produced by the bypass flow, $V_7$ is the bypass flow velocity at the outlet of the cold nozzle, $F_6$ is the hot thrust produced by the core flow, and $V_6$ is the velocity of the core flow at the exit of the core nozzle. $\dot{m}_0$, $V_0$, and $\dot{m}_{fuel}$, are ambient air intake mass flow rate, freestream velocity, and fuel mass flow rate, respectively.

### 3.2 Exergy based performance analysis

The second law of thermodynamics states that energy has quantity as well as quality. The quality of the energy is often measured in terms of exergy, which is the amount of energy that can be extracted as useful work [Hepbasli, 2008]. In general, the balance of exergy within a system in steady state can be expressed as follows:

$$\sum \chi_{in} - \sum \chi_{out} = \sum \chi_{des} \qquad (7)$$

where $\chi_{in}, \chi_{out}$ are the specific exergies flowing in and out of the system, while $\chi_{des}$ is the destruction of the exergy due to entropy generation in the system.

For a thermal system, the total exergy of a system is the sum of its potential, kinetic, physical, and chemical exergies [Turan et al., 2011]:

$$\chi = \chi_{pt} + \chi_{kn} + \chi_{ph} + \chi_{ch} \qquad (8)$$

The potential exergy is zero, as there is no distinctive elevation difference between the inlet and the outlet of the engine. The specific kinetic exergy is defined as follows:

$$\chi_{kn} = \frac{1}{2} V^2 \qquad (9)$$

The specific physical exergy is calculated as follows:

$$\chi_{ph} = h - h_0 - T_0 \left( \bar{s}_t - \bar{s}_0 - R \ln \frac{P}{P_0} \right) \qquad (10)$$

Since the entropy is the same at static and stagnation conditions, combining (8) and (9) yields:

$$\chi_{ph} + \chi_{kn} = h_t - h_{t,0} - T_{t,0} \left( \bar{s}_t - \bar{s}_{t,0} - R \ln \frac{P_t}{P_{t,0}} \right) \qquad (11)$$

where the subscript $t$ represents the stagnation condition of the properties, and $\bar{s}_t$ is the absolute entropy from (3) evaluated at its respective stagnation temperature.

The specific chemical exergies of liquid fuels can be obtained from the following expression [Canakci and Hosoz, 2006]:

$$\chi_{ch,fuel} = \left[1.0401 + 0.1728\frac{H}{C} + 0.0432\frac{O}{C} + 0.2169\frac{S}{C}\left(1 - 2.0628\frac{H}{C}\right)\right]|LHV| \qquad (12)$$

where $C$, $H$, and $S$ are the mass fractions of carbon, hydrogen, and sulphur, respectively. Liquid Kerosene with chemical equation $C_{12}H_{23.5}$ with LHV 43124kJ/kg has been assumed in this study. The air is assumed to consist of 77.48% $N_2$, 20.59% $O_2$, 0.03% $CO_2$, and 1.9% $H_2O$.

The combustion equation is therefore:

$$C_{12}H_{23.5} + 17.8750(O_2 + 3.7630 O_2 + 0.0015 CO_2 + 0.0923 H_2O) \rightarrow 12.0015 CO_2 + 62.2636 N_2 + 11.8423 H_2O \qquad (13)$$

The specific chemical exergy of the ideal gas mixture is:

$$\chi_{ch} = \frac{\bar{R}T_0}{M_f}\sum_{i=1}^{n} a_i \ln\left(\frac{y_i}{y_i^e}\right) \qquad (14)$$

where $\bar{R}$ = 8.314kJ/kgK is the universal gas constant, $M_f$ 167.8141kg/mole is the mass per mole of the Kerosene fuel, $y_i$ is the molar ratio of the $i^{th}$ component in the combustion products and $y_i^e$ is the molar ratio of the $i^{th}$ component in the reference environment (in this case the reference environment is the ambient environment). $a_i$ is the molar amount of the $i^{th}$ component. The resultant specific chemical exergy of the combustion product obtained from the above equations, as a function of $T_0$, is $4.5853T_0$ (J/kg).

The exergy efficiency, $\eta_{II}$ of the turbofan engine is evaluated with (15):

$$\eta_{II} = 1 - \sum_{i=1}^{n}\frac{\dot{m}_i \chi_{des,i}}{\dot{m}_{fuel}\chi_{ch,fuel}} \qquad (15)$$

4  Optimisation formulation

A computer program called TurboJet Engine Optimizer v1.0 (TJEO-1.0) was developed to perform the aforementioned analysis. The Pyevolve (an open source GA optimisation framework for Python programming language, available at http://pyevolve.sourceforge.net/) was used as the optimisation tool in TJEO-1.0. The optimisation problem is a single discipline, multi-variable multi-objective problem. It deals only with thermodynamic analysis.

There are eight design variables involved in the optimisation analysis as shown in Table 3. These variables are encoded into a chromosome-like string. Each of these parameters are represented by a group of 6-bit binary numbers within the string (also known as bit group). The precision of each of these groups is determined by the following equation [Homaifar et al., 1994]:

$$\Pi = \frac{U_{max} - U_{min}}{2^\lambda - 1} \qquad (16)$$

where $\Pi$ is the precision of the parameter, $U_{max}$ and $U_{min}$ are the upper and lower allowable values of the parameter, and $\lambda$ is the length of the group. Combining all the groups together, the resultant length of the chromosome string is 48-bit.

Thermodynamic constraints were taken into account in the optimisation. They are represented as hard constraints in GA. In other words, the fitness score of GA is penalised by setting its value to zero if one of the constraints is violated. The constraints that have been considered in this study are given in (17) and (18).

$$\pi_f \pi_{cL} \pi_{cH} < \pi_{max} \qquad (17)$$

$$P_{t,6} > P_0, \text{ and } P_{t,7} > P_0 \qquad (18)$$

Maximum stress on a turbofan engine happens during take-off, and this consideration is taken into account by setting the variables $H$ = 0, and $M$ = 0 and $T_{t,4}$ as 5% higher than the design point $T_{t,4}$. This setting is typical for contemporary jet engines [Borguet et al., 2006]. The values of other variables remained the same as those for cruise condition. The fitness score of GA is set to zero if any of the aforementioned constraints for take-off is violated.

Table 3. Design inputs and variables used in the simulation.

| Design Inputs | | Limits |
|---|---|---|
| $M$ | Freestream Mach number | ≤ 0.9 |
| $H$ | altitude | ≤ 12km |
| $\dot{m}_0$ | Inlet mass flow rate | ≥ 50kg/s |
| $\pi_{max}$ | Maximum compression ratio | ≤ 45 |
| Design Variables | | Constraints |
| $\alpha$ | Bypass ratio | 3.0 – 10.0 |
| $\beta$ | Bleed air fraction | 0.01 – 0.02 |
| $\varepsilon_1$ | Cooling air fraction for high-pressure turbine | 0.05 – 0.15 |
| $\varepsilon_2$ | Cooling air fraction for low-pressure turbine | 0.05 – 0.15 |
| $\pi_f$ | Fan pressure ratio | 1.2 – 2.0 |
| $\pi_{cL}$ | Compression ratio of LPC | 2.0 – 5.0 |
| $\pi_{cH}$ | Compression ratio of HPC | 4.0 – 10.0 |
| $T_{t,4}$ | Combustor outlet temperature | 1400 – 2000K |

As both the energy and exergy efficiencies used in this study are in efficiency terms, it is possible to formulate the multi-objective problem as a single objective problem with a scalar objective function. The objective function used to calculate the fitness score of solutions generated in each GA's iteration is formulated as follows:

$$score = \sqrt{\eta_I^2 + \eta_{II}^2} \qquad (19)$$

This function allows us to study the effects of each evaluation function $\eta_I$ and $\eta_{II}$, as well as the combination of both of them. When one of them is set to zero, the objective function is essentially the evaluation function which we wanted to study.

Roulette wheel selection has been chosen as the method to select chromosome in crossover within GA. The population size, generation size, and mutation probability are 160, 500, and 0.006, respectively.

## 5 Results and discussions

The design inputs selected for this study are: $H$ = 11km, $M$ = 0.86, $\dot{m}_0$ = 350kg/s, and $\pi_{max}$ = 45. Three sets of experiment were conducted (denoted as case (a), (b), and (c) hereafter) to evaluate how the exergy and energy efficiencies affect each other, as well as their effects on the engine performance. Case (a) was conducted with energy efficiency as the objective function, case (b) with exergy efficiency, and case (c) with the combination of energy and exergy efficiencies as the objective function. Each case was conducted with four repetitions. The results of 4 repeated GA runs (indicated by 4 line colours) are presented in Figure 2-9 where subfigures (a), (b), and (c) correspond to each of the aforementioned cases. The time taken for each run was 1060s on a 64-bit Operating system PC outfitted with a 2.67GHz CPU and 4GB of Random Access Memory (RAM).

- *Energy efficiency.* Figure 2 shows the energy efficiencies of the fittest individuals for each case. For case (a), the energy efficiency is the highest, around 40.5-41.0 % after 500 generations. For case (b), where exergy efficiency is used as objective function, the energy efficiency is the lowest, around 18.5-19.0%. This indicates that these two efficiencies contradict each other. Case (c) gives the optimum solutions to case (a) and case (b), with the energy efficiency hangs in between 38 and 40%, a small reduction in energy efficiency compared with case (a) and a big improvement by 200% compare with case (b).
- *Exergy efficiency.* In terms of exergy efficiency (as shown in Figure 3), case (a) has the lowest values among the three cases, around 62-64%. The exergy efficiency for case (b) is the highest, about 72.5%. This again confirms that the energy and exergy efficiencies contradict each other. For case (c), the exergy efficiency stays in

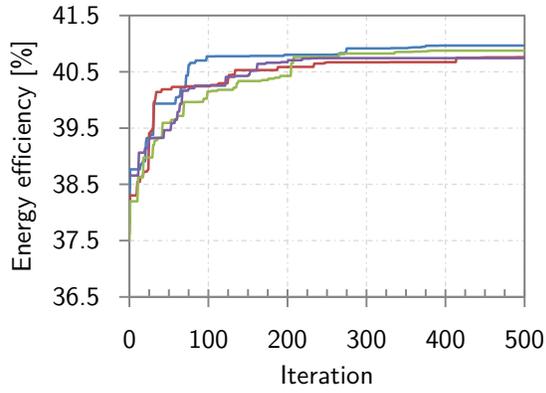

(a)

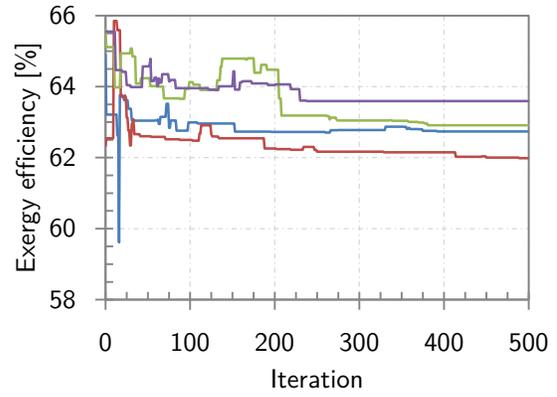

(a)

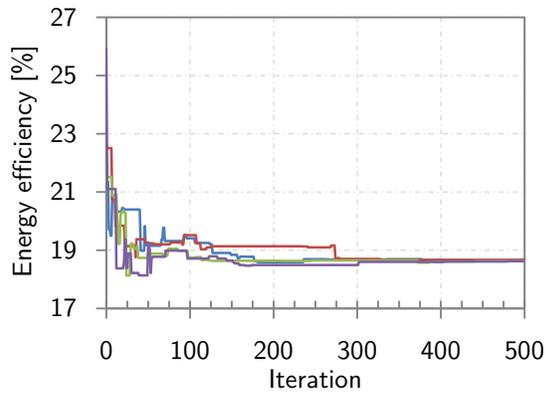

(b)

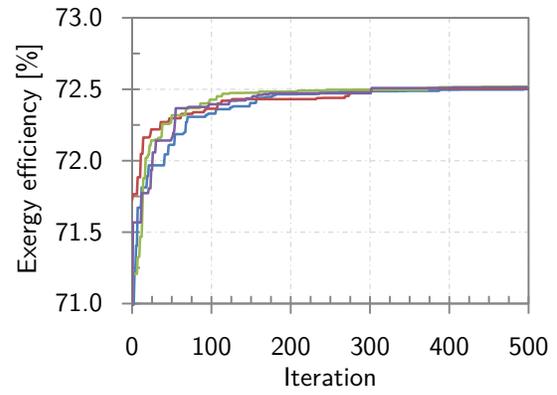

(b)

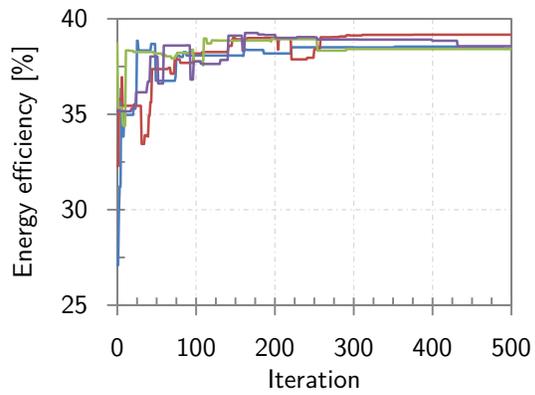

(c)

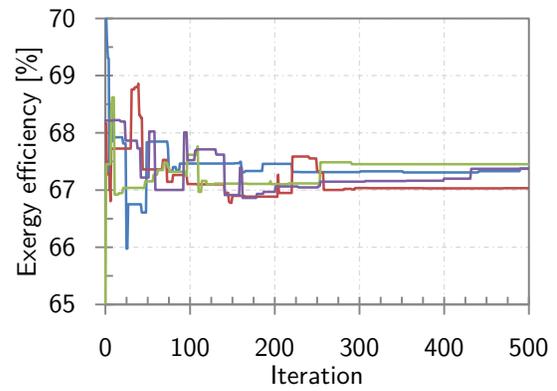

(c)

Figure 2. Convergence of energy efficiency: (a) $\eta_I$ as objective function, (b) $\eta_{II}$ as objective function, and (c) combination of $\eta_I$ and $\eta_{II}$ as objective function.

Figure 3. Convergence of exergy efficiency: (a) $\eta_I$ as objective function, (b) $\eta_{II}$ as objective function, and (c) combination of $\eta_I$ and $\eta_{II}$ as objective function.

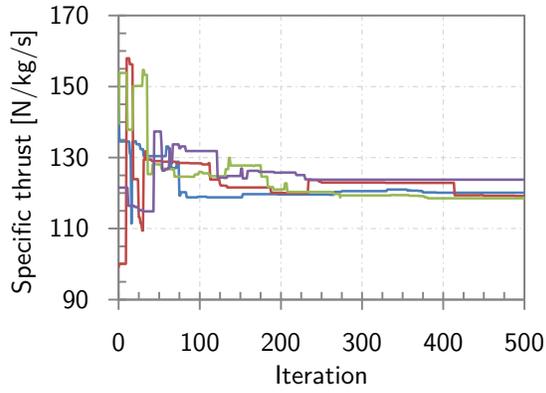

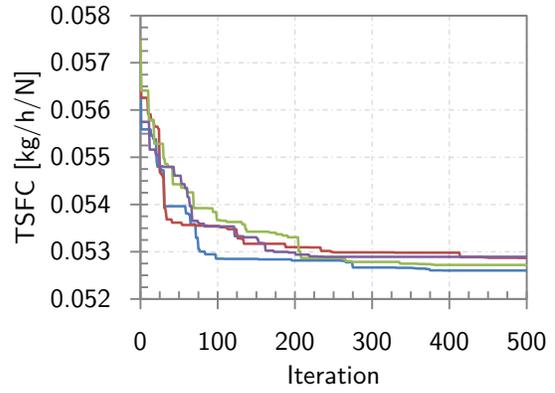

(a)

(a)

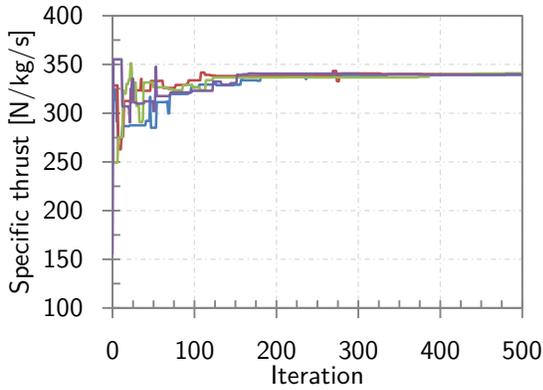

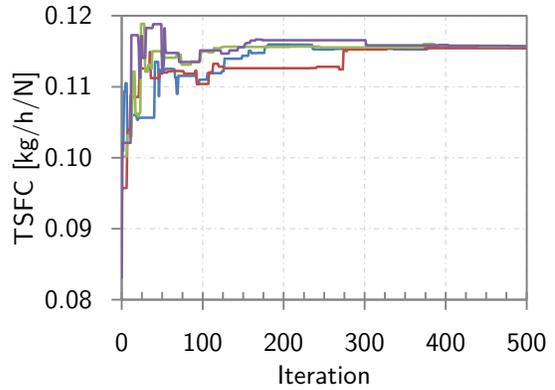

(b)

(b)

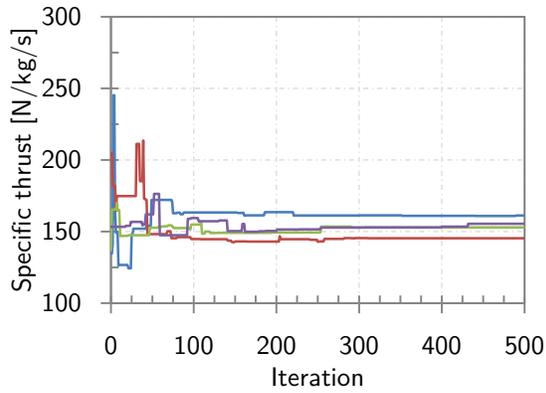

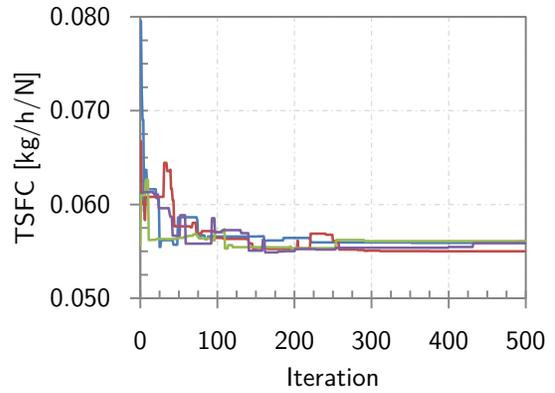

(c)

(c)

Figure 4. Convergence of specific thrust: (a) $\eta_I$ as objective function, (b) $\eta_{II}$ as objective function, and (c) combination of $\eta_I$ and $\eta_{II}$ as objective function.

Figure 5. Convergence of TSFC: (a) $\eta_I$ as objective function, (b) $\eta_{II}$ as objective function, and (c) combination of $\eta_I$ and $\eta_{II}$ as objective function.

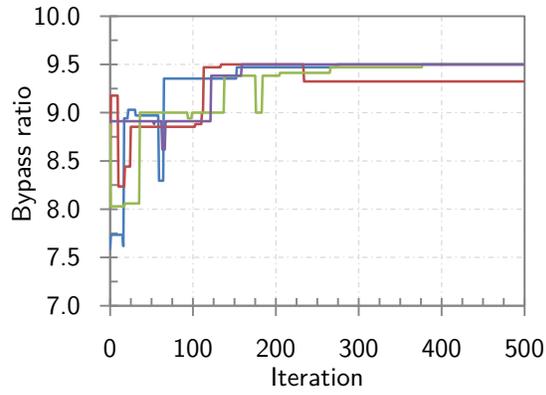

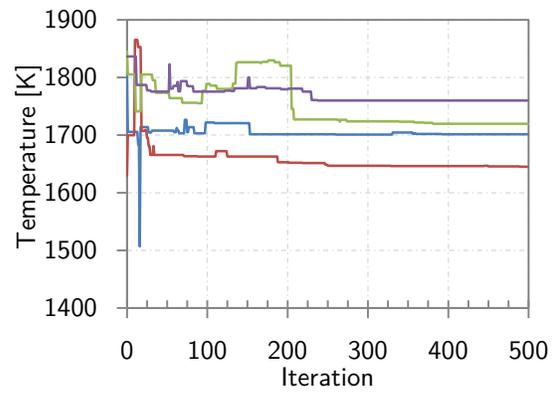

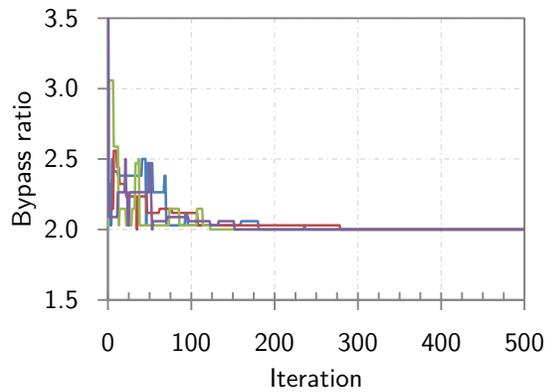

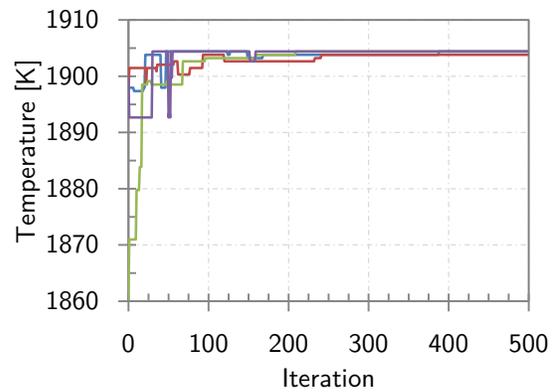

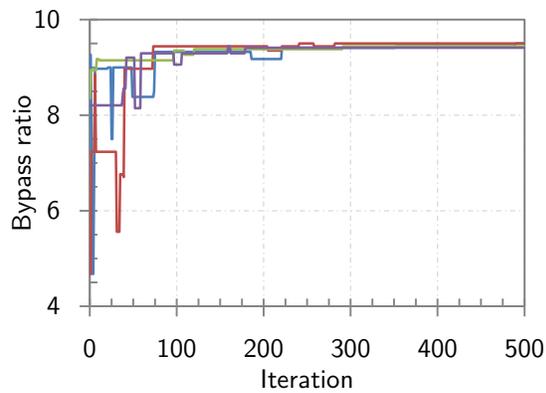

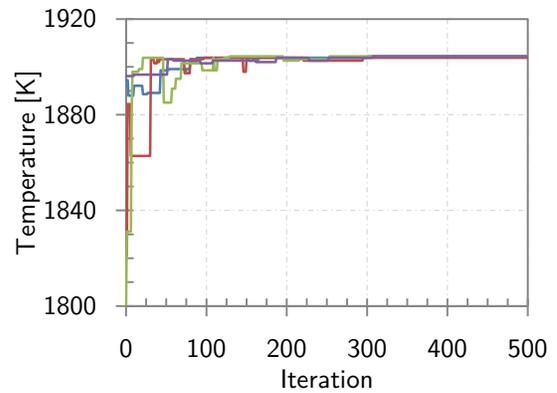

Figure 6. Convergence of bypass ratio: (a) $\eta_I$ as objective function, (b) $\eta_{II}$ as objective function, and (c) combination of $\eta_I$ and $\eta_{II}$ as objective function.

Figure 7. Convergence of burner outlet temperature $T_{t,4}$: (a) $\eta_I$ as objective function, (b) $\eta_{II}$ as objective function, and (c) combination of $\eta_I$ and $\eta_{II}$ as objective function.

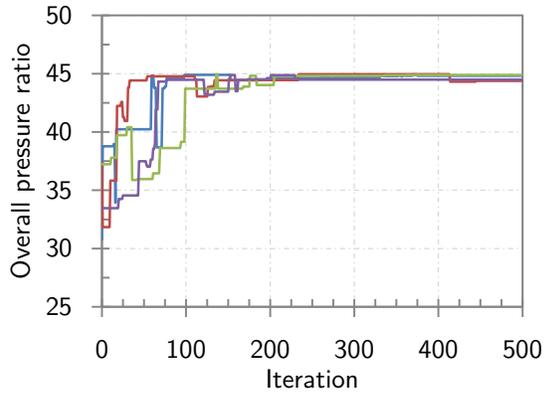

(a)

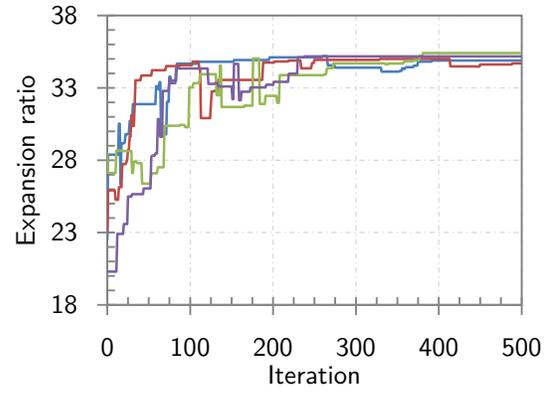

(a)

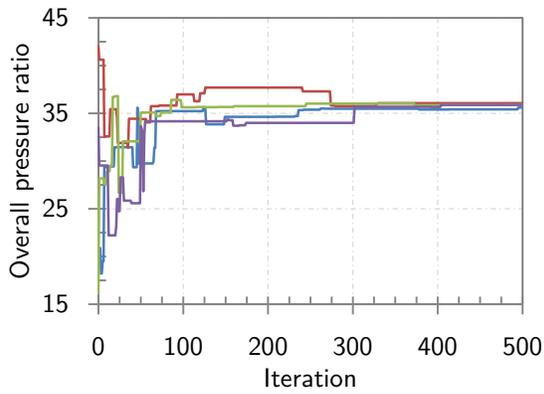

(b)

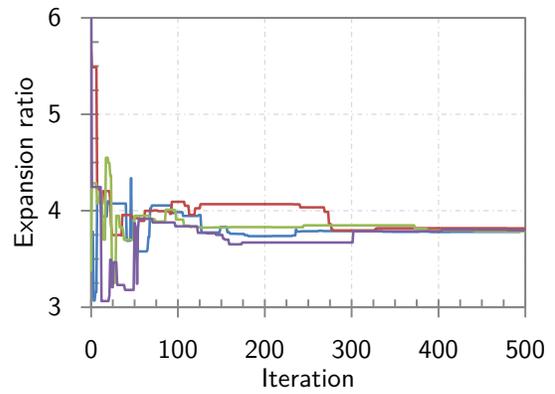

(b)

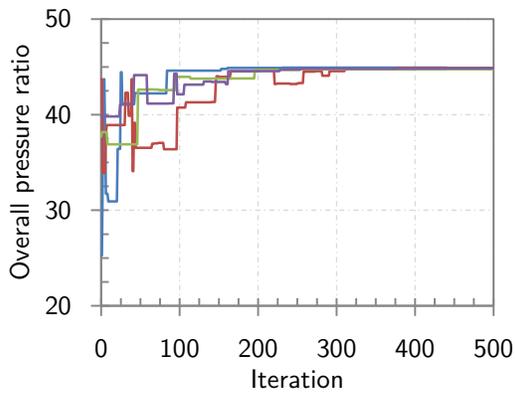

(c)

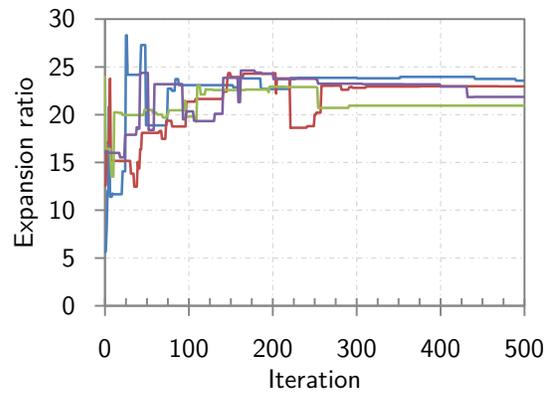

(c)

Figure 8. Convergence of overall pressure ratio: (a) $\eta_I$ as objective function, (b) $\eta_{II}$ as objective function, and (c) combination of $\eta_I$ and $\eta_{II}$ as objective function.

Figure 9. Convergence of turbine expansion ratio: (a) $\eta_I$ as objective function, (b) $\eta_{II}$ as objective function, and (c) combination of $\eta_I$ and $\eta_{II}$ as objective function.

between 67 and 68%, a solution which meets both the energy and exergy requirements.

- *Specific thrust.* From Figure 4, it appears that the optimum specific thrust for case (a) happens at around 118-125N/kg/s. On the other hand for case (b) when the exergy efficiency is chosen to be optimised, the specific thrust is concentrated at 340N/kg/s, meaning that the thrust produced per unit air intake is higher than that for case (a). As engine size is determined by the engine's intake area (which is determined by $\dot{m}_0$), higher specific thrust also means higher thrust to weight ratio can be achieved. For case (c), again the resultant specific thrust lies in between the values for case (a) and case (b), around 145-160N/kg/s.
- *Thrust specific fuel consumption.* Referring to Figure 5, the TSFC for case (a) is the lowest, about 0.0526-0.0528kg/h/N. For case (b), the value is the highest, about 0.115kg/h/N. For case (c), the resultant TSFC is rather low, about 0.055-0.056kg/h/N, slightly higher than that for case (a) but significantly lower than case (b). Although the TSFC for case (a) is the lowest, the specific thrust is not optimised. In case (c), 6% increase of TSFC yields an increase of specific thrust by 22-35%. This again justified the hypothesis the authors made earlier that exergy has to be taken into account for optimum design of jet engines. For case (c), the TSFC obtained is too high and not desirable.
- *Design variables*:
  - *Bypass ratio.* In terms of design parameters, the optimum $\alpha$ for case (a) and case (c) is 9.5, the highest possible value set in this study. This is consistent with the $\alpha$ of the contemporary turbofan engines for commercial jets, which are designed to minimise the fuel consumption. For case (b), $\alpha$ converges to 2.0, which is the lowest possible value for this study.
  - *Burner outlet temperature.* The $T_{t,4}$ for case (a) hangs in between 1640-1760K, the lowest among the three cases. For case (b) and case (c), the $T_{t,4}$ is invariably high, which is about 1905K (about 5% lower than the maximum $T_{t,4}$ set for this study). This infers that maximum $T_{t,4}$ has a role to play in exergy efficiency. With an increase in $T_{t,4}$, one can anticipate there will be an increase in the exergy efficiency.
  - *Overall pressure ratio.* The overall pressure ratio for case (a) and case (c) is the same, which is 45, the highest possible value used in this study. This infers that maximising overall pressure ratio tends to increase the energy efficiency of the engine. For case (b), the overall pressure ratio is 36, which is the lowest among the three cases.
  - *Turbine expansion ratio.* The turbine expansion ratio for case (a) is the highest among the three cases, around 35. For case (b), this value is the lowest, about nine times lower than case (a)'s. For case (c), it is around 20 to 25. From the authors' observation, energy efficiency increases with turbine expansion ratio. On the other hand, exergy efficiency decreases with turbine expansion ratio.

## 6 Conclusions

In this study, GA was successfully applied to a turbofan engine design problem, which was represented as a multi-objective and multi-variable optimisation problem. Energy and exergy efficiencies were used to formulate the objective functions for optimisation.

Three different sets of experiment had been carried out, each with $\eta_I$, $\eta_{II}$, and the combination of $\eta_I$ and $\eta_{II}$ as the objective function, respectively. The results showed that the combination of energy and exergy efficiencies as the objective function gave significant improvement on the exergy efficiency and specific thrust, although this accompanied by a small reduction in energy efficiency and a small increase in TSFC. This study also revealed that although minimum TSFC can be achieved by choosing the energy efficiency as the objective function, the work potential was not fully utilised.

Next, the authors suggest exergy efficiency to be included as one of the performance measurements in optimal engine design. This will help to maximise the work potential of the selected fuel.

Acknowledgements

The authors thank the Norwegian University of Science and Technology and the Brunel University for all support given in this research.

Turan O. and Karakoc T. H. Exergetic and energetic response surfaces for small turbojet engine. Applied Mechanics and Materials 110-116: 1054–1058, 2011.

Walsh P. P. and Fletcher P. Gas turbine performance, 2nd edition. Blackwell Science, 2004.

## Appendix A. Nomenclature

| Nomenclature | | | |
|---|---|---|---|
| $c_p$ | = | Specific heat at constant pressure | J/kgK |
| $F$ | = | Thrust | N |
| $f$ | = | Fuel-air ratio (fuel fraction) | |
| $H$ | = | Altitude | m |
| $h$ | = | Specific enthalpy | J/kg |
| $h_t$ | = | Total specific enthalpy | J/kg |
| $LHV$ | = | Lower heating value of fuel | kJ/kg |
| $M$ | = | Free stream Mach number | |
| $M_f$ | = | Mass per mole of fuel | kg/mole |
| $\dot{m}$ | = | Mass flow rate | kg/s |
| $P$ | = | Static pressure | Pa |
| $P_t$ | = | Total pressure (Stagnation pressure) | Pa |
| $R$ | = | Gas constant | J/kgK |
| $\bar{R}$ | = | Universal gas constant, 8.314 J/kgK | |
| $s$ | = | Specific entropy | J/kgK |
| $\bar{s}$ | = | Absolute specific entropy | J/kgK |
| $T$ | = | Static temperature | K |
| $T_t$ | = | Total temperature (Stagnation temperature) | K |
| $U$ | = | Value of a parameter in a bit group, see eq.(16) | |
| $V$ | = | Velocity | m/s |
| $W$ | = | Work | J |
| $\alpha$ | = | Bypass ratio | |
| $\beta$ | = | Bleed air fraction | |
| $\epsilon$ | = | Error | |
| $\varepsilon_1$ and $\varepsilon_2$ | = | Cooling air fraction for HPT and LPT, respectively | |
| $\gamma$ | = | Heat capacity ratio | |
| $\lambda$ | = | Length of a bit group, see eq.(16) | |
| $\eta$ | = | Component efficiency | |
| $\eta_I$ and $\eta_{II}$ | = | Engine's energy and exergy efficiencies, respectively | |
| $\Pi$ | = | Precision of a variable in a bit group, see eq.(16). | |
| $\pi$ | = | Pressure ratio | |
| $\chi$ | = | Specific exergy | J/kg |
| Subscript | | | |
| $b$ | = | burner | |
| $c$ | = | Core flow | |
| $cL$ | = | Low pressure compressor | |
| $cH$ | = | High pressure compressor | |
| $des$ | = | Destroyed | |

|       |   |                         |
|-------|---|-------------------------|
| $f$   | = | fan                     |
| $fuel$ | = | fuel                   |
| $i$   | = | $i^{th}$ component      |
| $in$  | = | Component / system inlet |
| $m$   | = | Temporary value for calculation |
| $max$ | = | Maximum value           |
| $min$ | = | Minimum value           |
| $out$ | = | Component / system outlet |
| *Superscript* |   |              |
| $*$   | = | Isentropic process      |

# Appendix B. Turbofan model

The following is the algorithm used to construct of the turbofan model which is used to calculate the properties of the working fluid at each engine station. This model is called by the GA to provide all the information needed for optimisation calculations.

*Compressor and fan:*
*Isentropic calculation:*
$\gamma = 1.4$
$T^*_{t,m} = T_{t,in}\, \pi^{(\gamma-1)/\gamma}$
$\epsilon = 1.0$
while ($\epsilon > 1 \times 10^{-5}$):
$\quad T_{t,m} = (T_{t,out} + T_{t,in})/2.0$
$\quad c_{p,m} = c_p(T_{t,m}, 0)$
$\quad \gamma = c_{p,m}/(c_{p,m} - R)$
$\quad T^*_{t,out} = T_{t,in}\, \pi^{(\gamma-1)/\gamma}$
$\quad \epsilon = |T^*_{t,m} - T^*_{t,out}|/T^*_{t,e}$
$\quad T^*_{t,m} = T^*_{t,out}$
$h^*_{t,out} = h(T^*_{t,out}, 0)$
$\Delta h^*_t = h^*_{t,out} - h_{t,in}$

*Polytropic calculation:*
$T_{t,out} = T^*_{t,out}$
$\epsilon = 1.0$
while ($\epsilon > 1 \times 10^{-5}$):
$\quad T_{t,m} = (T_{t,out} + T_{t,in})/2.0$
$\quad c_{p,m} = c_p(T_{t,m}, 0)$
$\quad \gamma = c_{p,m}/(c_{p,m} - R)$
$\quad \eta = (\pi^{(\gamma-1)/\gamma} - 1)/(\pi^{(\gamma-1)/e_c\gamma} - 1)$
$\quad T_{t,m} = \Delta h^*_t/(\eta\, c_{p,m}) + T_{t,in}$
$\quad \epsilon = |T_{t,out} - T_{t,m}|/T_{t,out}$
$\quad T_{t,out} = T_{t,m}$
$h_{t,out} = h(T_{t,out}, 0)$

*Burner:*

$f = 0.02$
$\epsilon = 1.0$
$while\ (\epsilon > 1 \times 10^{-5})$:
$\qquad h_{t,out} = h(T_{t,4}, f)$
$\qquad f_m = (h_{t,out} - h_{t,in})/(\eta_b |LHV|)$
$\qquad \epsilon = |f_m - f|/f$
$\qquad f = (f + f_m)/2$

*Turbine:*

$h_{t,out} = h_{t,out} - \sum \dot{W}_{out}/\dot{m}$
$T_{t,out} = T_{t,in}$
$\epsilon = 1.0$
$while\ (\epsilon > 1 \times 10^{-5})$:
$\qquad h_{t,m} = h(T_{t,out}, f_i)$
$\qquad \epsilon = |h_{t,m} - h_{t,out}|/h_{t,out}$
$\qquad c_{p,out} = c_p(T_{t,out}, f_i)$
$\qquad T_{t,out} = T_{t,out} + (h_{t,out} - h_{t,m})/c_{p,out}$